# Secrets and Quantifiers

**Béla Bajnok and Peter E. Francis**

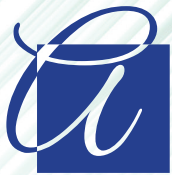

American writer, teacher, and comedian Sam Levenson gives us the following warning.

*Somewhere on this globe, every ten seconds, there is a woman giving birth to a child. She must be found and stopped.*

The joke comes from the fact that the statements "for any given time, there is a place, at which there is a woman giving birth" and "there is a place, for which there is a woman that is always giving birth" are not equivalent. The first is believable while the second, if true, would be ridiculously disastrous.

As this example shows, there is a discrepancy between the sometimes deceiving and confusing use of quantifiers (words or phrases used to describe the quantity of objects that have a certain property) in the English language and the rigidly fundamental role that quantification plays in mathematical writing. So the question is: how do we translate?

## Some Notations

We write $P(x)$ to denote a *predicate*, which is a statement whose truth depends on $x$ in a set $A$. For example, if $P(x)$ denotes the predicate "$x < 57$ or $x$ is odd," for a positive integer $x$, then $P(4)$, $P(57)$, and $P(101)$ are all true, while $P(100)$ and $P(134)$ are both false.

The symbol $\forall$ represents the universal quantifier (pronounced "for all"). So the sentence

$$\forall x \in A, P(x)$$

means "for all $x$ in the set $A$, $P(x)$ is true." The symbol $\exists$ is the existential quantifier (pronounced "there exists"). The sentence

$$\exists x \in A, P(x)$$

means "there is some (at least one) $x$ in the set $A$ for which $P(x)$ is true."

The order and type of nested quantifiers can change the meaning of the statement drastically. Sam Levenson's joke plays off the differences between the two statements

$\forall$ time, $\exists$ place, $\exists$ woman giving birth, and
$\exists$ place, $\exists$ woman, $\forall$ time giving birth.

## A Secret Game

Imagine that Suzy has a secret sequence of four positive integers that Quentin would like to know. He asks questions in the form of sequences of positive integers, to which Suzy's response will be the scalar product (the sum of the coordinate-wise products) of the two sequences: if Suzy's secret sequence is $\mathbf{s} = (s_1, s_2, s_3, s_4)$, the response to the question $\mathbf{q} = (q_1, q_2, q_3, q_4)$ is

$$\mathbf{q} \cdot \mathbf{s} = q_1 s_1 + q_2 s_2 + q_3 s_3 + q_4 s_4.$$

Quentin has an idea of how to determine Suzy's sequence. By brute force, he can ask the four questions (2, 1, 1, 1), (1, 2, 1, 1), (1, 1, 2, 1), and (1, 1, 1, 2). Suzy responds with $w$, $x$, $y$, and $z$, respectively, and these answers determine the following system of four linear equations in the variables $s_1$, $s_2$, $s_3$, and $s_4$.

$$\begin{cases} w = 2s_1 + s_2 + s_3 + s_4 \\ x = s_1 + 2s_2 + s_3 + s_4 \\ y = s_1 + s_2 + 2s_3 + s_4 \\ z = s_1 + s_2 + s_3 + 2s_4 \end{cases}$$

Maybe Quentin has studied some linear algebra or maybe he was lucky, because these four questions are *linearly independent*: the system of equations has a unique solution given by

$$\begin{cases} s_1 = (4w - x - y - z)/5 \\ s_2 = (4x - w - y - z)/5 \\ s_3 = (4y - w - x - z)/5 \\ s_4 = (4z - w - x - y)/5. \end{cases}$$

It is exciting (at least for Quentin) that with four questions, he can always discover the secret. However, you might be wondering if Quentin can do better. Are fewer questions sufficient? This possibility takes a bit more to unpack.

We can think of some secrets where single questions decode them. For example, if Quentin asks the question (1, 5, 10, 20) and Suzy answers 36, the secret must be (1, 1, 1, 1): there is only one way to use $1, $5, $10, and $20 bills to pay $36 so that each currency is used at least once. The situation is the same if Suzy answers 37, 38, 39, or 40; however, if Suzy answers 41, her secret could be (6, 1, 1, 1) or (1, 2, 1, 1). To go further, we will need to be more technical with what we mean by 'decoding the secret with one question.'

### Decoding Sequences

Let $D(\mathbf{q}, \mathbf{s})$ denote the predicate "the question $\mathbf{q}$ decodes the secret $\mathbf{s}$." To decode the secret with one question, we must ensure that no other secret sequence returns the same response for the question. That is, $D(\mathbf{q}, \mathbf{s})$ is true exactly when there is no sequence $\mathbf{t}$ different from $\mathbf{s}$ for which $\mathbf{q} \cdot \mathbf{s} = \mathbf{q} \cdot \mathbf{t}$. Equivalently in quantifier notation,

$$\forall \mathbf{t},\ \mathbf{t} \neq \mathbf{s},\ \mathbf{q} \cdot \mathbf{t} \neq \mathbf{q} \cdot \mathbf{s}.$$

But is a single question always enough? The answer depends largely on the two possible precise interpretations of this question:
1) Is "$\exists \mathbf{q},\ \forall \mathbf{s},\ D(\mathbf{q},\mathbf{s})$" true?
2) Is "$\forall \mathbf{s},\ \exists \mathbf{q},\ D(\mathbf{q},\mathbf{s})$" true?

Our examination of these two questions in the next sections exposes that the order of the quantifiers results in two rather different outcomes.

### The Master Key

The first quantifier chain "$\exists \mathbf{q},\ \forall \mathbf{s},\ D(\mathbf{q},\mathbf{s})$" means "there is some fixed question that can decode any secret." Imagine questions as keys and secrets as locks; the predicate $D(\mathbf{q}, \mathbf{s})$ can be interpreted as "key $\mathbf{q}$ opens lock $\mathbf{s}$."

Using the metaphor, the previous statement means "there is a key that can open any lock." If this were true in the game context, Quentin would always win by asking the same one question, making the game quite boring. Is this true?

Let $\mathbf{q} = (q_1, q_2, q_3, q_4)$ be an arbitrary question. Define $\mathbf{s} = (1, 1, 1+q_4, 1)$ and $\mathbf{t} = (1, 1, 1, 1+q_3)$. Then

$$\mathbf{q} \cdot \mathbf{s} = q_1 + q_2 + q_3(1+q_4) + q_4$$
$$= q_1 + q_2 + q_3 + q_4(1+q_3)$$
$$= \mathbf{q} \cdot \mathbf{t}.$$

Thus, given a question, we can build at least two sequences between which the question cannot distinguish. In other words, if we thought we had a "master key" (one that could open any lock), we would be wrong because we showed how to build at least two locks that the key cannot open.

### The Unbreakable Secret

The second quantifier chain "$\forall \mathbf{s},\ \exists \mathbf{q},\ D(\mathbf{q},\mathbf{s})$" translates to "for any secret, we can pick a question to decode it." This means there is no secret safe from being decoded in one question. While this sentence sounds quite similar to the previous, "master key," sentence, the small change in quantifiers can make a big difference in meaning. We showed that the statement in 1) is false, but the statement in 2) is true. Let's see why.

Suppose that $\mathbf{s} = (s_1, s_2, s_3, s_4)$ is arbitrary. Pick four pairwise relatively prime positive integers $a_1, a_2, a_3,$ and $a_4$ greater than $s_1, s_2, s_3,$ and $s_4$, respectively. Then let

$$q_1 = a_2 a_3 a_4,\ q_2 = a_1 a_3 a_4,\ q_3 = a_1 a_2 a_4,\ q_4 = a_1 a_2 a_3,$$

and $\mathbf{q} = (q_1, q_2, q_3, q_4)$. We claim $\mathbf{q}$ decodes $\mathbf{s}$. To show this, assume there is some sequence $\mathbf{t} = (t_1, t_2, t_3, t_4)$ with $\mathbf{q} \cdot \mathbf{s} = \mathbf{q} \cdot \mathbf{t}$. Then we have

$$0 = q_1(s_1 - t_1) + q_2(s_2 - t_2) + q_3(s_3 - t_3) + q_4(s_4 - t_4).$$

Because $a_1$ divides the last three summands and 0, $a_1$ must also divide $q_1(s_1 - t_1)$. As $a_1$ is relatively prime to $q_1$, $a_1$ must divide $s_1 - t_1$. Both $s_1$ and $t_1$ are positive integers, so $s_1 - t_1 < s_1$. Because $a_1 > s_1$, $a_1$ can only divide $s_1 - t_1$ if $t_1 \geq s_1$. Similarly, $t_2 \geq s_2$, $t_3 \geq s_3$, and $t_4 \geq s_4$. As the previous equation shows four nonpositive integers summing to 0, we conclude that each summand must be 0, so $s_i = t_i$ for each $i$. Thus, $\mathbf{s} = \mathbf{t}$. There is no 'unbreakable secret'!

### The Solution

As we have seen, every secret can be unlocked by a single question (though we can't know the question in advance), and there are four questions that can unlock any secret. But can we use fewer questions? As before, the answer depends on how we make this precise, and we now have a new consideration: are the questions being asked at the same time, or are we aware of the answer to the first before asking the second? It turns out that any first question along with a carefully chosen follow-up are always sufficient. Can you prove it?

Hopefully you can see how these examples show the effect that quantifier order has on the truth of statements, as well as the change in caliber of proof needed to keep up with seemingly small variations in sentence structure. At the very least, now if anyone tells you that they have a shirt for every day of the week, you can rightfully ask them how they manage to wash it! ●

*Béla Bajnok is a professor of mathematics at Gettysburg College and is the director of the American Mathematics Competitions of the MAA.*

*Peter E. Francis is an undergraduate mathematics student at Gettysburg College.*